# Non-integrality of integrable over Kähler differential forms

Achim Hennings [♦]

**Abstract:** We show that for a normal isolated singularity the module of square-integrable regular differential forms of top degree is not integral over the module of Kähler differential forms. This is related to the fibre dimension of the Nash transform.

In [KW] it is shown that for a normal Cohen-Macaulay singularity of dimension $d \geq 2$ the canonical module consisting of (on the smooth locus) regular $d$−differential forms is not integral over the module of Kählerian differential forms.[1] In [SiSmU] is pointed out the relation to the fibre dimension of the Nash transformation, and an analogous assertion is proved for homogeneous singularities, but the general case is left open. Here we prove the same assertion for normal isolated singularities by showing that already the submodule of square-integrable $d$−forms is not integral over the module of Kähler forms.

Let $(X, 0) \subseteq (\mathbb{C}^N, 0)$ be a reduced irreducible and non-smooth isolated singularity of dimension $d \geq 2$. Let $X$ be a closed analytic representative with smooth part $X^* = X - \{0\}$ in a sufficiently small open ball $U \subseteq \mathbb{C}^N$.

Let $\nu: \bar{X} \to X$ be the normalization and $\pi: \tilde{X} \to \bar{X}$ a resolution of singularities. We have maps

$$\nu^* \Omega_X^d \to \pi_* \omega_{\tilde{X}} \to \omega_{\bar{X}}.$$

The kernel of the first map is the torsion while the second one is injective. The module $\omega_{\bar{X}}$ can be embedded as an ideal as follows. Let $f = (f_1, \ldots, f_c)$ with $c = N - d$ be a generic system of generators for the ideal of $X$ (i.e. elements of the ideal which form a sequence of parameters and generate the ideal except at an analytic subset not containing $X$). A meromorphic $d$−form with divisor $K \leq 0$ is given by

$$\alpha = \frac{dx_1 \ldots dx_N}{df_1 \ldots df_c} | X^*,$$

where $x_1, \ldots, x_N$ are coordinates. There are ideals $I \subseteq J \subseteq \mathcal{K} \subseteq \mathcal{O}_{\bar{X}}$ with the property $\omega_{\bar{X}} = \mathcal{K}\alpha$, $\pi_* \omega_{\tilde{X}} = J\alpha$ and $\text{image}(\nu^* \Omega_X^d) = I\alpha$. By means of

$$\omega_{\bar{X}} \xrightarrow{df_1 \ldots df_c} \nu^*(\Omega_U^N \otimes \mathcal{O}_X) \cong \mathcal{O}_{\bar{X}},$$

the form $\alpha$ is mapped to 1. So this map has image $\mathcal{K}$. We remark that the embedding of $\omega_{\bar{X}}$ is unique up to a rational function. Therefore we may define the notion of integral dependence among the modules considered here by the corresponding ideals.

**Theorem:** $\nu^* \Omega_X^d \to \pi_* \omega_{\tilde{X}}$ is not integral, i.e. $I \subseteq J$ is not integral.

---

[♦] Universität Siegen, Fakultät IV, Hölderlinstraße 3, D-57068 Siegen

[1] In fact, they show that the Kähler differentials are mapped into the product of the canonical module by the maximal ideal.

Proof: As the map factorizes over $\Omega_{\bar{X}}^d$, we may immediately pass to the normalization, as long as it is singular again. On the other hand, if $\bar{X}$ is smooth, the image of $\nu^*\Omega_X^d \to \omega_{\bar{X}} \cong \mathcal{O}_{\bar{X}}$ is the Jacobian ideal of the composition $\bar{X} \xrightarrow{\nu} X \to \mathbb{C}^N$ (in local coordinates) which is a proper ideal because $X$ is singular and irreducible.

Proceeding inductively, we consider intersections $X_0 := X \cap (x_1)$ by a general hyperplane. Then $X_0$ is again a true singularity with $X_0^* := X^* \cap X_0$ smooth. Furthermore, it is reduced and if $d \geq 3$ (as we may suppose for the induction) it is irreducible (cf. [BG], [F]). Let $\nu_0: \bar{X}_0 \to X_0$ be its normalization.

We have a commutative diagram:

$$\begin{array}{ccccc}
dx_1 \wedge \Omega_X^{d-1} \subseteq \Omega_X^d & \to & \omega_X & \xrightarrow{df_1\ldots df_c} & \Omega_U^N \otimes \mathcal{O}_X \cong \mathcal{O}_X \\
\downarrow & & /dx_1 \downarrow & & \downarrow \\
\nu_{0*}\nu_0^*\Omega_{X_0}^{d-1}/torsion & \to & \nu_{0*}\omega_{\bar{X}_0} & \xrightarrow{dx_1 df_1\ldots df_c} & \nu_{0*}\nu_0^*(\Omega_U^N \otimes \mathcal{O}_{X_0}) \cong \nu_{0*}\mathcal{O}_{\bar{X}_0}
\end{array}$$

The map on the right is the composition $\mathcal{O}_X \to \mathcal{O}_{X_0} \to \nu_{0*}\mathcal{O}_{\bar{X}_0}$ of the restriction with the normalization. In the first line the images of $dx_1 \wedge \Omega_X^{d-1} \subseteq \Omega_X^d$ in $\mathcal{O}_X$ are the determinantal ideals of size $c$

$$I_c(\partial f/\partial(x_2, \ldots, x_N))\mathcal{O}_X \subseteq I_c(\partial f/\partial(x_1, \ldots, x_N))\mathcal{O}_X.$$

In the second line the corresponding image is

$$I_c(\partial f/\partial(x_2, \ldots, x_N))\mathcal{O}_{\bar{X}_0}.$$

As the coordinates may be chosen generally, $\partial f/\partial x_1$ is integral on $X$ over the remaining partials $\partial f/\partial x_i$. To see this, we remark that the Rees ring of the module of partials has fibre dimension at most $N-1$, so this module has a reduction generated by $N-1$ elements. From the integrality of the module of all partial derivatives over the same module with $\partial f/\partial x_1$ left out, it can be deduced that the corresponding determinantal ideals of size $c$ are also integral (in $\mathcal{O}_X$). Therefore, the above two ideals have the same integral closure. (Cf. [Ga], [T])

By this observation, we have only to show, that $dx_1 \wedge \Omega_X^{d-1} \to \pi_*\omega_{\bar{X}}$ is not integral.

Let $\pi_0: X_0' \to \bar{X}_0$ be a resolution of singularities. Below, we choose the strict transform of $X_0$. We obtain an enlarged diagram:

$$\begin{array}{ccccc}
dx_1 \wedge \Omega_X^{d-1} & \xrightarrow{\alpha} & \pi_*\omega_{\bar{X}} & \to & \omega_X \subseteq \mathcal{O}_X \\
\gamma \downarrow & & \delta \searrow & & \downarrow \\
\nu_{0*}\nu_0^*\Omega_{X_0}^{d-1}/torsion & \xrightarrow{\beta} & \nu_{0*}\pi_{0*}\omega_{X_0'} & \to & \nu_{0*}\omega_{\bar{X}_0} \subseteq \nu_{0*}\mathcal{O}_{\bar{X}_0}
\end{array}$$

We claim the following assertions, which will be proved later.

<u>Assertion 1</u>: image $(\delta) \supseteq \nu_{0*}\pi_{0*}\omega_{X'_0}$.

<u>Assertion 2</u>: The theorem holds for $d = 2$.

For the proof of the theorem, we may assume by induction on the dimension with base given by assertion 2, that $\beta$ is not integral. All the more, $\nu_{0*}\nu_0^*\Omega_{X_0}^{d-1}/torsion \to$ image $(\delta)\nu_{0*}\mathcal{O}_{\tilde{X}_0}$ is not integral by assertion 1. Therefore, $\alpha$ cannot be integral.

It remains to prove the auxiliary assertions.

Proof of assertion 1: We choose the resolution in such a way that the preimage of $X_0$ is a divisor with normal crossings of smooth components:

$$\tilde{X}_0 = \pi^* X_0 = X'_0 + \sum_i e_i E_i$$

By multiplication with $t = \pi^* x_1$ we obtain the exact sequence

$$0 \to \omega_{\tilde{X}} \xrightarrow{t} \omega_{\tilde{X}} \to \omega_{\tilde{X}_0} \to 0.$$

By adjunction we have $\omega_{\tilde{X}_0} = \omega_{\tilde{X}}(\tilde{X}_0)/\omega_{\tilde{X}}$ and $\omega_{X'_0} = \omega_{\tilde{X}}(X'_0)/\omega_{\tilde{X}}$, so $\omega_{X'_0} \subseteq \omega_{\tilde{X}_0}$. The theorem of Grauert-Riemenschneider [GR] states that $H^1(\tilde{X}, \omega_{\tilde{X}}) = 0$. Accordingly $\Gamma(\omega_{\tilde{X}}) \to \Gamma(\omega_{\tilde{X}_0})$ is surjective and $\Gamma(\omega_{X'_0}) \subseteq \Gamma(\omega_{\tilde{X}_0})$, in particular forms in $\Gamma(\omega_{X'_0})$ are extendable.

Proof of assertion 2: We can take $X$ to be normal. Let $\pi: \tilde{X} \to X$ be the minimal resolution and $E = \bigcup E_i$ the exceptional set. For $\Omega_X^2 \to \pi_*\omega_{\tilde{X}}$ to be not integral it is sufficient that (the pullback of) all forms in $\Omega_X^2$ vanish along $E$, but some sections of $\omega_{\tilde{X}}$ do not.

Assertion 3: Let $p \in E$ and $\varphi \in \Omega_{X,0}^2$. Then $\pi^*\varphi \in \omega_{\tilde{X},p}$ has a zero.

Namely the family of all $\pi^*\varphi$ generates, expressed in local coordinates about $p$, in $\omega_{\tilde{X},p} \cong \mathcal{O}_{\tilde{X},p}$ the Jacobian ideal of the composition $\tilde{X} \to X \to \mathbb{C}^N$, which has a zero because otherwise this would be an immersion at $p$.

Assertion 4: If $M = \tilde{X}$ is the minimal resolution, there is an exceptional component $E_i$ and a form $\varphi \in \Gamma(M, \omega_M)$ with $\operatorname{ord}_{E_i}\varphi = 0$.

Proof: This is a direct consequence of the preparatory considerations in [L] which we recall in the following. The canonical bundle on $M$ is denoted here by $K$. We have $K \cdot E_i \geq 0$ for all components because $M$ is minimal, and this is essentially used.

Let $Z$ be the uniquely determined fundamental cycle, i.e. $Z > 0$ is minimal satisfying $Z \cdot E_i \leq 0$ for all $E_i$. Let $Z_0, Z_1, \ldots, Z_l = Z$ be a computation sequence,

$$Z_0 = 0, Z_1 = Z_0 + E_{i_1}, Z_k = Z_{k-1} + E_{i_k}, Z_{k-1} \cdot E_{i_k} > 0, 2 \leq k \leq l.$$

As $Z$ is unique and $Z^2 < 0$, $E_{i_1}$ can be chosen such that $Z \cdot E_{i_1} < 0$. One defines new sequences as

$$E_{i_1} = E_{j_l}, \ldots, E_{i_l} = E_{j_1},$$

$$Y_0 = Z - Z_l = 0, Y_1 = Z - Z_{l-1} = E_{j_1}, \ldots, Y_l = Z - Z_0 = E_{j_1} + \cdots + E_{j_l}.$$

Let $E_i$ be an exceptional component and $\nu: \bar{E}_i \to E_i$ its normalization. Let $L$ be a line bundle on $E_i$. Let

$$g' = h^1(\mathcal{O}_{E_i}) = g + \delta, g = h^1(\mathcal{O}_{\bar{E}_i}), \delta = h^0(\nu_*\mathcal{O}_{\bar{E}_i}/\mathcal{O}_{E_i}).$$

By Riemann-Roch we have:

Theorem 1 [L, (2.2)]: If $c(L) > 2g' - 2$, then $h^1(\mathcal{O}(L)) = 0$ and $h^0(\mathcal{O}(L)) = c(L) + 1 - g'$.

For a line bundle $L$ on $M$ one can deduce:

Theorem 2 [L, prop. 2.1]: If for all exceptional components $L \cdot E_i \geq K \cdot E_i$ and if $Z \cdot E_{i_1} < 0$ by the above choice, then for $n \geq 0, 1 \leq k \leq l$

$$H^1(\mathcal{O}(L - nZ - Y_{k-1})) = 0.$$

The first theorem can be applied for $L = K|E_i$, for all $i$. In fact $g' = 1 + \frac{1}{2}(K \cdot E_i + E_i^2)$ and therefore $c(K|E_i) = K \cdot E_i > 2g' - 2 = K \cdot E_i + E_i^2$. One gets $h^1(\mathcal{O}_{E_i}(K)) = 0$ and

$$h^0(\mathcal{O}_{E_i}(K)) = K \cdot E_i + 1 - g' = \frac{1}{2}(K \cdot E_i - E_i^2) > 0 \text{ as } K \cdot E_i \geq 0.$$

The second theorem we apply for $L = K$ and $nZ + Y_{k-1} = E_{j_1}$ (i.e. $n = 0$, $k = 2$ if $l > 1$ and $n = 1$, $k = 1$ if $l = 1$), and obtain from $H^1(\mathcal{O}(K - E_{j_1})) = 0$ the exact sequence with non-vanishing right term

$$0 \to H^0(\mathcal{O}(K - E_{j_1})) \to H^0(\mathcal{O}(K)) \to H^0\left(\mathcal{O}_{E_{j_1}}(K)\right) \to 0.$$

This concludes the proof of assertion 4 and 2.